\documentclass[journal,usenames,dvipsnames]{IEEEtran}

\usepackage{myStyleIEEE}
\usepackage{bm}
\usepackage{smartdiagram}
\usesmartdiagramlibrary{additions}
\usepackage{steinmetz}
\usepackage{eurosym}
\begin{document}
	
\title{Active Distribution Grids offering Ancillary Services in Islanded and Grid-connected Mode}

\author{
\IEEEauthorblockN{Stavros~Karagiannopoulos,~\IEEEmembership{Student Member,~IEEE}, Jannick Gallmann, Marina Gonz\'alez Vay\'a, Petros~Aristidou,~\IEEEmembership{Member,~IEEE,} and~Gabriela~Hug,~\IEEEmembership{Senior Member,~IEEE}}%
  \thanks{S. Karagiannopoulos and G. Hug are with the Power Systems Laboratory, ETH Zurich, 8092 Zurich, Switzerland. Email: \{karagiannopoulos $|$ hug\}@eeh.ee.ethz.ch.}%
  \thanks{P. Aristidou is with the School of Electronic and Electrical Engineering, University of Leeds, Leeds LS2 9JT, UK. Email: p.aristidou@leeds.ac.uk}%
  \thanks{J. Gallmann is with the Stadtwerk Winterthur, Switzerland, Email: jannick.gallmann@alumni.ethz.ch}
  \thanks{M. Gonz\'alez Vay\'a is with EKZ, Switzerland. Email: Marina.GonzalezVaya@ekz.ch}
}


\maketitle

\IEEEpeerreviewmaketitle
\begin{abstract}
Future active distribution grids (ADGs) will incorporate a plethora of Distributed Generators (DGs) and other Distributed Energy Resources (DERs), allowing them to provide ancillary services in grid-connected mode and, if necessary, operate in an islanded mode to increase reliability and resilience. In this paper, we investigate the ability of an ADG to provide frequency control (FC) in grid-connected mode and ensure reliable islanded operation for a pre-specified time period. First, we formulate the operation of the grid participating in European-type FC markets as a centralized multi-period optimal power flow problem with a rolling horizon of 24 hours. Then, we include constraints to the grid-connected operational problem to guarantee the ability to switch to islanded operation at every time instant. Finally, we explore the technical and economic feasibility of offering these services on a balanced low-voltage distribution network. \review{The results show that the proposed scheme is able to offer and respond to different FC products, while ensuring that there is adequate energy capacity at every time step to satisfy critical load in the islanded mode.}
\end{abstract}


\section{Introduction}
While moving towards a low-carbon, sustainable electricity system, future Distribution Networks (DNs) are expected to host a large share of Distributed Generators (DGs) to satisfy the demand currently supplied by fossil-fuel and nuclear power plants. DGs, coordinated with other Distributed Energy Resources (DERs), such as electric vehicles, Battery Energy Storage Systems (BESSs) and Flexible Loads (FLs), consequently amplify the role of DNs, making them an important part in ensuring grid reliability and resilience~\cite{Panteli2017}, and enabling them to provide ancillary services to transmission voltage levels~\cite{PES-TR22}. Thus, it is crucial to operate modern DNs \emph{actively}, i.e. controlling DERs to ensure secure, reliable and cost-effective operation. 

In this paper, we consider a centralized method with existing communication infrastructure, which is a valid assumption in modern DNs that do not cover large geographical areas~\cite{olivares2014trends}. 
\vspace{-0.2cm}
\subsection{Related work}
Operating active DNs using optimization has been widely explored in literature, e.g.~\cite{tsikalakis2011centralized,mitra2012determination,Falahi2013dynamic,khodaei2014microgrid,moya2008distributed,asano2009optimization,joos2000potential,yuen2007feasibility,rueda2013capacity,koller2015review, koller2016primary,Majzoobi2017}. Here, we only review work concerned with the DN capability to operate off-grid and the provision of ancillary services offered in grid-connected mode.
\subsubsection{Islanded operation}\label{subsec:isl_op}
In~\cite{mitra2012determination}, a Monte Carlo (MC) approach is applied in the design stage to determine the required BESS size to reliably operate in this mode but without incorporating the possibility of the BESS offering ancillary services in the connected mode. Reference~\cite{khodaei2014microgrid} on the other hand optimizes the microgrid operational costs in grid-connected mode as a master problem, while ensuring islanded capability for multiple hours as a subproblem. However, neither the provision of energy based ancillary services nor the incorporation of voltage control, which requires the consideration of a network model, are addressed. A model-predictive-control \review{(MPC)} scheme, including the dynamics of the system, is used in \cite{Falahi2013dynamic} to predict future voltage instabilities and adjust the reactive power generation accordingly. Here, the focus is only on keeping voltages close to nominal values in islanded mode, not offering other grid-connected or islanded services. Further,~\cite{koller2015review} examines the behavior of a real BESS offering frequency control reserves and supporting islanded operation. The described setup uses a dispatchable diesel generator in addition to intermittent renewable energy sources but no network constraints are considered and the response of the BESS is based on heuristics rather than on centralized optimization.

\subsubsection{Grid-connected operation}\label{subsec:pw_AS}
In grid-connected mode, the main objective is usually to operate the DN in the most cost-effective way. \review{A detailed review of the state-of-the-art research in microgrids is presented in~\cite{Parhizi2015}, where the authors review around $400$ works, covering the areas of
microgrid economics, operation, control, protection, and
communications. Reference~\cite{tsikalakis2011centralized} investigates the economic evaluation of grid-connected microgrids that participate in real-time markets but without considering islanded operation. Ref.~\cite{Majzoobi2017} focuses also on the optimal scheduling of an active DN providing frequency regulation, load leveling and ramping services. Offering ancillary services by various DER technologies is explored in~\cite{moya2008distributed}, while~\cite{asano2009optimization,joos2000potential} analyze the economic feasibility and the potential amount of reserve provision by distributed generation. However, they do not consider BESSs or include network constraints. The technical feasibility of providing ancillary services with multiple microgrids as a pool bidder is investigated in~\cite{yuen2007feasibility}, while~\cite{rueda2013capacity} takes the perspective of a Transmission System Operator (TSO), minimizing its own expenditures and evaluating the competitiveness of DGs. 
Finally, several DSOs are already providing actual frequency control products to the TSO. For instance~\cite{koller2015review} and~\cite{koller2016primary} discuss the operational experience of a BESS offering primary frequency control in the European interconnected network. However neither of~\cite{asano2009optimization,joos2000potential,Majzoobi2017,koller2015review,koller2016primary} considers the network modeling within the optimization.}

\review{The consideration of grid-connected and islanded mode of an active distribution grid in combination with the provision of ancillary service taking into account also grid constraints is, to the best of our knowledge, not considered in any previous work.}

\vspace{-0.2cm}
\subsection{Contributions}
In this paper, we propose a centralized optimization approach to operate an Active Distribution Grid (ADG). We explicitly incorporate uncertainty into the formulation and consider the opportunity of offering ancillary services. The proposed methodology is based on a multi-period, Chance-Constrained Optimal Power Flow (CC-OPF) formulation, where various frequency control products are offered by DERs. We ensure that at any point the DN can operate in islanded mode for a limited time by including additional constraints in the centralized problem. In this way, the DN operation considers both the uncertainties from RES as well as the potential need for islanded operation~\cite{JannickCIRED}.

Consequently, the main contributions of this paper are:
\begin{itemize}
    \item A multi-period CC-OPF formulation that:
        \begin{itemize}
            \item allows an ADG to offer ancillary services in grid-connected mode, while being able to switch to islanded mode at any time, and
            \item considers RES uncertainty through a rolling horizon strategy.
        \end{itemize}
    \item A case-study analyzing various frequency control (FC) products and the performance of the proposed method in offering these products.
\end{itemize}

\review{The remainder of the paper is organized as follows: Section~\ref{OPFall} presents the mathematical formulation of the deterministic OPF considering ancillary service provision as well as islanded operation, while Section~\ref{Sec:uncertainty_finalproblem} accounts for uncertainty and presents the final CC-OPF} Then, Section~\ref{case} introduces the case study and simulation results for the islanded and grid-connected case. Finally, conclusions are drawn in Section~\ref{Conclusion}.
\vspace{-0.2cm}
\review{\section{Centralized deterministic OPF}\label{OPFall}
In this section, we present the deterministic centralized OPF scheme used to compute the optimal DER setpoints. The objective considers both the grid-connected and the islanded mode simultaneously and is optimized in a rolling horizon fashion; we model FC products offered in the grid-connected mode as constraints, while at the same time enabling a potential switch to islanded mode for the following 24 hours.}

\subsection{Centralized OPF}\label{detOPF}
\review{\subsubsection{Preliminaries}
We consider a radial balanced distribution grid with $\mathcal{N}$ being the set of nodes \review{using the index $j$}, $\mathcal{T}$ the set of branches \review{using the index $i$}, $\mathcal{B} \subseteq \mathcal{N}$ the subset of nodes with BESS, $\mathcal{L} \subseteq \mathcal{N}$ the subset of nodes with loads, $\mathcal{F} \subseteq \mathcal{L} \subseteq \mathcal{N}$ the subset of flexible (controllable) loads, and $\mathcal{R} \subseteq  \mathcal{N}$ the subset of nodes with DGs. The DER control measures (detailed below) are represented by the variable $u$, and the variables referring to the islanded operation mode by the subscript ``$\textrm{isl}$''.}
\subsubsection{Objective function}\label{sec:obj} The objective function is defined as
\begin{equation}\label{eq:obj}
\min \limits_{\mathbf{u}}  \sum_{t=t_{\textrm{MPC}}}^{t_{\textrm{MPC}} + T } \left(  C_{\textrm{t}}^{\textrm{curt,g}} + C_{\textrm{t}}^{\textrm{curt,l}} + C_{\textrm{t}}^{\textrm{exc}} + C_{\textrm{t}}^{\textrm{AS}} \right)  \Delta t,
\end{equation}
\noindent where $t_{\textrm{MPC}}$ denotes the current time step of the MPC algorithm, $T$ the rolling horizon period and $\Delta t$ the length of a time interval within the horizon.

At each discrete time $t$, the objective function consists of four terms;
\paragraph{\text{\normalfont\ $C_{\textrm{t}}^{\textrm{curt,g}}$}} This term corresponds to the cost of generation curtailment in both the grid-connected and islanded mode and is given by
\begin{equation}
    C_{\textrm{t}}^{\textrm{curt,g}} = \sum_{j \in \mathcal{R}} c_{\textrm{t}}^{\textrm{curt,g}} \cdot ( P_{\textrm{j,t}}^{\textrm{curt,g}} + f_{\textrm{isl}} \cdot P_{\textrm{j,t}}^{\textrm{curt,g,isl}}), 
\end{equation}
\noindent where $P_{\textrm{j,t}}^{\textrm{curt,g}}=P_{\textrm{j,t}}^{\textrm{g,max}}-P_{\textrm{j,t}}^{\textrm{g}}$ is the curtailed power of the DG connected at node $j$ at time $t$ (resp., $P_{\textrm{j,t}}^{\textrm{curt,g,isl}}$ in the islanded case), $P_{\textrm{j,t}}^{\textrm{g,max}}$ the maximum available active power, and $P_{\textrm{j,t}}^{\textrm{g}}$ the actual in-feed; $c_{\textrm{t}}^{\textrm{curt,g}}$ is the cost of curtailment at time $t$, and $f_{\textrm{isl}}$  a constant scalar that adjusts the cost in the islanded case. \review{The cost of generation curtailment is policy-related in the grid-connected case and can be very different from country to country. Typically, generators are compensated at the prevailing electricity market price, whereas in some European countries, they are compensated only for a small part of the curtailed energy~\cite{Lew2013}. In California, compensation for curtailment begins after a contractually agreed number of hours which vary among contracts~\cite{Bird2014}. However, in the islanded case, the operation of the DGs becomes more important, since they are the only sources to satisfy the local demand, i.e. no external grid is available. Thus, there is another value associated with the injection of power from DGs in the islanded case, which is accounted for by the scalar $f_{\textrm{isl}}$.} 

\paragraph{\text{\normalfont\ $C_{\textrm{t}}^{\textrm{curt,l}}$}} This term represents the cost of load curtailment in the islanded mode and is given as
\review{\begin{equation}
C_{\textrm{t}}^{\textrm{curt,l}} = \sum_{j \in \mathcal{L}} c_{\textrm{t}}^{\textrm{curt,l,isl}} \cdot P_{\textrm{j,t}}^{\textrm{curt,l,isl}}, 
\end{equation}
\noindent where $c_{\textrm{t}}^{\textrm{curt,l,isl}}$ is the cost of load curtailment at time $t$, and $P_{\textrm{j,t}}^{\textrm{curt,l,isl}}$ the curtailed load. In the grid-connected case, we assume that any local generation-load mismatch can be covered from the transmission network without the need for load shedding in the grid-connected case, similar to~\cite{Romero-Quete2018}}.

\paragraph{\text{\normalfont\ $C_{\textrm{t}}^{\textrm{exc}}$}} The third term includes the cost of exchanging power with the upper voltage levels and is given by
\begin{equation}
C_{\textrm{t}}^{\textrm{exc}} = c_{\textrm{t}}^{\textrm{buy}} \cdot P_{\textrm{1,t}}^{\textrm{buy}} - c_{\textrm{t}}^{\textrm{sell}} \cdot P_{\textrm{1,t}}^{\textrm{sell}}, 
\end{equation}
where $c_{\textrm{t}}^{\textrm{buy}}$ ($c_{\textrm{t}}^{\textrm{sell}}$) is the price of buying (selling) electric energy from (to) the main grid. By considering different prices for buying and selling, \review{i.e. at each time step buying electricity is more expensive than selling,}  we prioritize storing excess energy locally (promoting the self-consumption of the DN), over exporting power to higher voltage levels; $P_{\textrm{1,t}}^{\textrm{g}}=P_{\textrm{1,t}}^{\textrm{buy}} - P_{\textrm{1,t}}^{\textrm{sell}}$ ($P_{\textrm{1,t}}^{\textrm{buy}},P_{\textrm{1,t}}^{\textrm{sell}}\geq0$) is the active power exchange measured at the substation \review{making sure that the ADG cannot buy and sell electricity at the same time. A similar formulation is followed in~\cite{stavrosEEM} to determine the position, i.e. short or long, of an aggregator participating in energy markets}. 
\paragraph{\text{\normalfont\ $C_{\textrm{t}}^{\textrm{AS}}$}} The final term corresponds to revenues from offering ancillary services to upper voltage levels, given by
\begin{equation}
C_{\textrm{t}}^{\textrm{AS}} = c_{\textrm{t}}^{\textrm{bid}} \cdot P_{\textrm{bid}}, 
\end{equation}
\noindent where $c_{\textrm{t}}^{\textrm{bid}}$ is the pay-as bid volume-weighted average price of the accepted bids in the frequency control market from the respective week of the previous year (assumed known) and $P_{\textrm{bid}}$ the bid (to be determined by the CC-OPF).

\subsubsection{Power balance constraints}\label{powerbalance}
 The power injections at every node $j$ and time step $t$ are given by
 \begin{subequations}
\label{eq:nodal_active_power_balance}
\begin{gather}
P_{\textrm{j,t}}^{\textrm{inj}}= P_{\textrm{j,t}}^{\textrm{g}} - P_{\textrm{j,t}}^{\textrm{lflex}}  - \left( P_{\textrm{j,t}}^{\textrm{B,ch}} - P_{\textrm{j,t}}^{\textrm{B,dis}} \right), 
\\
\label{eq:nodal_reactive_power_balance}
Q_{\textrm{j,t}}^{\textrm{inj}}= Q_{\textrm{j,t}}^{\textrm{g}} + Q_{\textrm{j,t}}^{\textrm{B}} - P_{\textrm{j,t}}^{\textrm{lflex}} \cdot \tan(\phi_{\textrm{l}}),
\end{gather}
\end{subequations}
where $P_{\textrm{j,t}}^{\textrm{g}}$ and $Q_{\textrm{j,t}}^{\textrm{g}}$ are the active and reactive power injections of the generators at node $j$; $P_{\textrm{j,t}}^{\textrm{lflex}}$ and $P_{\textrm{j,t}}^{\textrm{lflex}} \cdot \tan(\phi_{\textrm{l}})$ are the active and reactive node demands (after control), with $cos(\phi_{\textrm{l}})$ being the power factor of the load; $Q_{\textrm{j,t}}^{\textrm{B}}$ is the reactive power of the BESS and, $P_{\textrm{j,t}}^{\textrm{B,ch}}$ and $P_{\textrm{j,t}}^{\textrm{B,dis}}$ are respectively the charging and discharging BESS active powers.

\subsubsection{Power flow constraints}\label{powerflow}
In this work, we integrate the Backward/Forward Sweep (BFS) \review{method} into our \review{power flow} formulation~\cite{stavrosPowertech,StavrosIREP,StavrosPSCC18}. \review{The solution of the BFS \emph{power flow} problem is achieved  iteratively, by "sweeping" the distribution network and updating the network variables at each iteration, which consists of two sweeps. First, in the backward sweep step, the current injections at all buses and the corresponding branch currents are calculated. Then, in the forward sweep step, the currents are used to calculate the voltage drop over all branches, updating the bus voltages for the next iteration of the algorithm. Within an OPF framework, we consider only one iteration to model network flows and to avoid the non-linearities introduced by the AC power flow equations. Subsequently, if the derived solution is not AC feasible, we update the voltages by projecting the solution into the AC feasible manifold~\cite{StavrosIREP}, and re-run the OPF problem.} This reformulation provides a sufficiently accurate approximation of the full AC OPF~\cite{Fortenbacher2016a}, is computationally tractable~\cite{StavrosPSCC18}, and results in AC feasible solutions which can account for uncertainties (see~\cite{StavrosIREP} for more details). A single iteration of the BFS is used to replace the AC power-flow constraints in the OPF formulation as follows:
\begin{subequations}\label{eq:}
\begin{gather}
    I_\textrm{j,t}^{\textrm{inj}}= \left (\frac{(P_{\textrm{j,t}}^{\textrm{inj}} + jQ_{\textrm{j,t}}^{\textrm{inj}})^{*}}{\bar{V}_{\textrm{j,t}}^{*}}\right),   \\
    I_\textrm{t}^{\textrm{br}}=BIBC \cdot I_\textrm{t}^{\textrm{inj}},  \\
    \Delta{V}_{\textrm{t}}=BCBV \cdot I_\textrm{t}^{\textrm{br}},  \\
    V_{{\textrm{j,t}}}= V_\textrm{slack} - \Delta V_\textrm{tap} \cdot \rho_{\textrm{t}} + \Delta{V}_{\textrm{t}},  \\
     \rho_{min} \leq \rho_{\textrm{t}} \leq \rho_{max}, \label{eq:OLTC3}  
\end{gather}
\end{subequations}
\noindent where $\bar{V}_{\textrm{j,t}}^{*}$ is the voltage magnitude at node $j$ at time $t$,$~^*$~indicates the complex conjugate and the bar indicates that the value from the previous iteration is used; $I_\textrm{t}^{\textrm{inj}}=[ I_\textrm{j,t}^{\textrm{inj}}, j\hspace{-0.1cm}\in\hspace{-0.1cm}\mathcal{N} ]$ and $I_\textrm{t}^{\textrm{br}}=[ I_\textrm{i,t}^{\textrm{br}}, i\hspace{-0.1cm}\in\hspace{-0.1cm}\mathcal{T} ]$ represent respectively the vectors of bus injection and branch flow currents; $I_\textrm{i,t}^{\textrm{br}}$ is the $i$-th branch current; $BIBC$ is a matrix with ones and zeros, capturing the radial topology of the DN; the entries in $\Delta{V}_{\textrm{t}}$ correspond to the voltage drops over all branches and phases; $BCBV$ is a matrix with the complex impedances of the lines as elements; $V_\textrm{slack}$ is the voltage in per unit at the slack bus (here assumed to be \review{$1\phase{0^{\circ}}$}); $\Delta V_{tap}$ is the voltage magnitude change caused by one tap action of the On-Load Tap Changer (OLTC) transformer and assumed constant for all taps for simplicity; and, $\rho_{\textrm{t}}$ is an integer value defining the position of the OLTC position. The parameters ($\rho_\textrm{min},\rho_\textrm{max}$) are respectively the minimum and maximum tap positions of the OLTC transformer.

\subsubsection{Thermal loading and voltage constraints} \label{VIconstraints}
The constraints for the current magnitudes at time $t$ are given by
\begin{align}
    |I_\textrm{i,t}^{\textrm{br}}|   \leq I_{\textrm{i}}^{\textrm{max}}, 
\end{align}
\noindent where $I_{\textrm{i}}^{\textrm{max}}$ is the maximum thermal limit for the $i$-th branch.

Similarly, the voltage constraints are given by
    $$V_\textrm{min}   \leq | V_{\textrm{j,t}} | \leq V_{\textrm{max}}. $$ 
\noindent where $(V_{\textrm{max}}$, $V_{\textrm{min}})$ are respectively the upper and lower acceptable voltage limits. However, the lower voltage magnitude limit results in a non-convex constraint~\cite{StavrosPSCC18}. By exploiting the fact that the voltage angles are typically small in distribution grids, we can approximate the complex voltage with its real part for the lower bound, as explained in~\cite{StavrosPSCC18}:
\begin{align}
    V_\textrm{min}  & \leq {Re}\left\{ V_{\textrm{j,t}}   \right\} , \qquad  | V_{\textrm{j,t}} | \leq V_{\textrm{max}},  
\end{align}

\subsubsection{DER constraints}\label{DERineq}
\paragraph{DG limits}

In this work, we consider inverter-based DGs such as PVs. Their limits are given by
\begin{align}
    P_{\textrm{j,t}}^{\textrm{g,min}} \leq P_{\textrm{j,t}}^{\textrm{g}} \leq P_{\textrm{j,t}}^{\textrm{g,max}}, \quad Q_{\textrm{j,t}}^{\textrm{g,min}} \leq  Q_{\textrm{j,t}}^{\textrm{g}} \leq Q_{\textrm{j,t}}^{\textrm{g,max}}, 
\end{align}
\noindent where $P_{\textrm{j,t}}^{\textrm{g,min}}$, $P_{\textrm{j,t}}^{\textrm{g,max}}$, $Q_{\textrm{j,t}}^{\textrm{g,min}}$ and $Q_{\textrm{j,t}}^{\textrm{g,max}}$ are the upper and lower limits for active and reactive DG power at each node $j \in \mathcal{N}$ and time $t$. These limits vary depending on the type of the DG and the control schemes implemented. Usually, small DGs have technical or regulatory~\cite{VDE} limitations on the power factor they can operate at or reactive power they can produce. Any of these limitations can be captured in this constraint.  

\paragraph{Controllable loads} 
We further consider flexible loads which can shift a limited amount of energy consumption in time. The loads are therefore modeled by
\begin{subequations}\label{eq:CL}
\begin{gather}
    P_{\textrm{j,t}}^{\textrm{lflex}} = P_{\textrm{j,t}}^{\textrm{l}} + f_{\textrm{j,t}}^{\textrm{lflex}} \cdot  P_{\textrm{j,t}}^{\textrm{shift}}, \\
   -1 \leq f_{\textrm{j,t}}^{\textrm{lflex}} \leq 1,\\
   \sum \limits_{t=t_{\textrm{0}}}^{t_{\textrm{MPC}}-1} f_{\textrm{j,t}}^{\textrm{lflex}} + \sum \limits_{t=t_{\textrm{MPC}}}^{t_{\textrm{MPC}}+ T} f_{\textrm{j,t}}^{\textrm{lflex}} =0, \label{CLsum}
\end{gather}
\end{subequations}
\noindent where 
$P_{\textrm{j,t}}^{\textrm{shift}}$ is the shiftable load of the non-shiftable demand $P_{\textrm{j,t}}^{\textrm{l}}$; $f_{\textrm{j,t}}^{\textrm{lflex}}$ is the normalized factor defining the final load shift. The past values for $f_{\textrm{j,t}}^{\textrm{lflex}}$, i.e. for time instances $t=t_{\textrm{0}}$ (start of the simulation) to $t_{\textrm{MPC}}-1$, are constant. This is necessary due to the moving horizon approach and the fact that the total load at the end of the simulation period needs to be maintained which is ensured by \eqref{CLsum}. \review{The separation of these terms is done for clarity reasons, to distinguish the fixed past values from the decision variables of the optimization problem.} 

\paragraph{Battery Energy Storage Systems}
Finally, the constraints related to the BESS at node $j$ are given as
\begin{subequations}
\label{eq:BESS}
\begin{gather}
    SoC_{\textrm{min}}^{\textrm{B}} \cdot E_{\textrm{cap,j}}^{\textrm{B}} \leq  E_{\textrm{j,t}}^{\textrm{B}} \leq SoC_{\textrm{max}}^{\textrm{B}} \cdot E_{\textrm{cap,j}}^{\textrm{B}},      \label{eq:BESS_en}\\  
    E_{\textrm{j,1}}^{\textrm{B}} = E_{\textrm{j,$t_{\textrm{0}}$}},    \label{eq:BESS_SOC}\\
    E_{\textrm{j,t}}^{\textrm{B}}  = E_{\textrm{j,t-1}}^{\textrm{B}} + (\eta_{\textrm{B}} \cdot  P_{\textrm{j,t}}^{\textrm{B,ch}} - \frac{P_{\textrm{j,t}}^{\textrm{B,dis}}}{\eta_{\textrm{B}}}) \cdot \Delta t ,   \label{eq:dynBESS1}\\
    0 \leq P_{\textrm{j,t}}^{\textrm{B,ch}} \leq P_{\textrm{j,max}}^{\textrm{B}} ,  \quad  0 \leq P_{\textrm{j,t}}^{\textrm{B,dis}} \leq P_{\textrm{j,max}}^{\textrm{B}},  \label{eq:BESS_P} \\
    P_{\textrm{j,t}}^{\textrm{B,ch}} +   P_{\textrm{j,t}}^{\textrm{B,dis}} \leq  \textrm{max}(P_{\textrm{j,t}}^{\textrm{B,ch}}, P_{\textrm{j,t}}^{\textrm{B,dis}}),\label{eq:BESS_P2} \\
    |Q_{\textrm{j,t}}^{\textrm{B}}| \leq \max{ \left\{ P_{\textrm{j,t}}^{\textrm{B,ch}},P_{\textrm{j,t}}^{\textrm{B,dis}} \right\}} \cdot tan(\phi_{\textrm{max}}^{\textrm{B}}), \label{eq:BESS_Qb}
\end{gather}
\end{subequations}
where $E_{\textrm{cap,j}}^{\textrm{B}}$ is the installed BESS capacity connected at node $j$; 
$SoC_{\textrm{min}}^{\textrm{B}}$ and $SoC_{\textrm{max}}^{\textrm{B}}$ are the fixed minimum and maximum per unit limits for the battery state of charge; and, $E_{\textrm{j,t}}^{\textrm{B}}$ is the available energy at node $j$ and time $t$. The initial energy content of the BESS in the first time period is given by $E_{\textrm{j,$t_{\textrm{0}}$}}$, and (\ref{eq:dynBESS1}) updates the energy in the storage at each period $t$ based on the BESS efficiency $\eta_{\textrm{B}}$, time interval $\Delta t$ and the charging and discharging power of the BESS $P_{\textrm{j,t}}^{\textrm{B,ch}}$ and $ P_{\textrm{j,t}}^{\textrm{B,dis}}$. The charging and discharging powers are defined as positive according to \eqref{eq:BESS_P}, while \eqref{eq:BESS_P2} ensures that the BESS is not charging and discharging at the same time. Finally, \eqref{eq:BESS_Qb} limits the reactive power output as a function of the charging or discharging power and the maximum power
factor $cos(\phi_{\textrm{max}}^{\textrm{B}})$;
\vspace{-0.4cm}
\subsection{Ancillary services}\label{subsec:AS}
In grid connected mode, we include the offering of frequency control products following a European market framework~\cite{ENTSO-E2009, swissgrid}. These require power and energy reserves, that can be called at any time. In the following sections, we describe the technical constraints of each product. Please note that only one single FC product is offered at a time, i.e. multiple services are not considered.

\subsubsection{Primary frequency control (PFC)}\label{subsubsec:FCR}
\review{PFC is a symmetrical product, i.e. each bid needs to provide symmetrical power bands both for up- and down-regulation, to cover imbalances both from excess production or consumption. }
The European frequency control reserve cooperation~\cite{ENTSO-E2009} has set the energy requirement to $0.25 \cdot P_{\textrm{bid}}$, i.e. the provider has to be able to deliver the full committed power ($P_{\textrm{bid}}$) for a quarter of an hour ($15$ minutes). However, evaluation of realized primary control signals showed that this requirement is conservative\cite{JannickMA}, i.e. much less energy is actually needed. In this work, only the battery is considered to be able to offer this product. 
 The power reserves for up- and down-regulation are given by 
\begin{subequations}\label{eq:FCR_power_down}
\begin{gather}
\sum_{j\in\mathcal{B}} \left( P_{\textrm{max,j}}^{\textrm{B}} - P_{\textrm{j,t}}^{\textrm{B,dis}} + P_{\textrm{j,t}}^{\textrm{B,ch}}  \right) \geq P_{\textrm{bid}},  \\
\sum_{j\in\mathcal{B}} \left( P_{\textrm{max,j}}^{\textrm{B}} - P_{\textrm{j,t}}^{\textrm{B,ch}} + P_{\textrm{j,t}}^{\textrm{B,dis}} \right) \geq P_{\textrm{bid}},
\end{gather}
\end{subequations}
\noindent where $P_{\textrm{bid}}$ is the weekly power size of the PFC bid. 

The energy that has to be reserved $\forall t$ is given by
\begin{subequations}\label{eq:FCR_en_up}
\begin{gather}
\sum_{j\in\mathcal{B}} \left( E_{\textrm{j,t}}^{\textrm{B}} - SoC_{\textrm{min}}^{\textrm{B}} \cdot E_{\textrm{cap,j}}^{\textrm{B}}\right) \geq P_{\textrm{bid}} \cdot \Delta t_{\textrm{1}}, \\
\sum_{j\in\mathcal{B}} \left( SoC_{\textrm{max}}^{\textrm{B}} \cdot E_{\textrm{cap,j}}^{\textrm{B}} - E_{\textrm{j,t}}^{\textrm{B}} \right) \geq P_{\textrm{bid}} \cdot \Delta t_{\textrm{1}},
\end{gather}
\end{subequations}
\noindent where $\Delta t_{\textrm{1}}$ is defined to be 15 minutes \cite{Swissgrid_AS_principles}.

\subsubsection{Secondary Frequency Control (SFC)}\label{subsubsec:SFC}
\review{SFC is activated after PFC to bring frequency back to the nominal value, and restore the scheduled power exchanges with other control areas. SFC is also symmetrical and requires fast response times. Thus,} for the provision of this product, we employ the BESS and the PV units. The power reserves for up- and down-regulation (again symmetrical) are given by
\begin{subequations}\label{eq:SFC_power_up}
\begin{gather}
\sum_{j\in\mathcal{B}} \left( P_{\textrm{max,j}}^{\textrm{B}} - P_{\textrm{j,t}}^{\textrm{B,dis}} + P_{\textrm{j,t}}^{\textrm{B,ch}} \right) + \sum_{j\in\mathcal{R}} \left( P_{\textrm{j,t}}^{\textrm{g,max}} - P_{\textrm{j,t}}^{\textrm{g}} \right) \geq P_{\textrm{bid}}, \\
\label{eq:SFC_power_down}
\sum_{j\in\mathcal{B}} \left( P_{\textrm{max,j}}^{\textrm{B}} - P_{\textrm{j,t}}^{\textrm{B,ch}} + P_{\textrm{j,t}}^{\textrm{B,dis}} \right) + \sum_{j\in\mathcal{R}} \left( P_{\textrm{j,t}}^{\textrm{g}} \right) \geq P_{\textrm{bid}},
\end{gather}
\end{subequations}
\noindent where $P_{\textrm{bid}}$ is the weekly power size of the SFC bid. Secondary control is activated after a few seconds and is typically completed after 15 minutes~\cite{swissgrid}. \review{However, in reality this scheme does not guarantee that the energy requirement will not exceed the energy required for a provision of $P_{\textrm{bid}}$ for 15 minutes. By design, there is a continuous secondary call signal that needs to be followed, not accounting for specific energy requirements. For this reason, a statistical approach was followed to analyze ex-post the SFC signal over 1 year in Switzerland, and subsequently derive hourly worst case energy requirements per bid size of secondary frequency power. }
The worst case values for a 24-h rolling horizon required an energy content of around 5.5 hours times the amount of the bid size in either direction~\cite{JannickMA}. However, these values are too conservative, and would limit drastically the flexibility on the secondary frequency control market. Thus, we consider as additional constraints only the first $4$ hours of the worst case requirements. Afterwards, the missing/surplus energy can still be bought/sold at the spot market with a lead time of one hour~\cite{EPEX}.

Furthermore, since PV forecasts are subject to some short-term adjustments, we require that at least 50\% of the energy of a worst case call has to come from the BESS.  

Thus, the energy content evolution for the first 4 hours of the worst case call is described by
\begin{multline}\label{eq:SFC_en_up}
E_{\textrm{j,t+$\vartheta$}}^{\textrm{B,2,+}} = \underbrace{E_{\textrm{j,t+$\vartheta$-1}}^{\textrm{B,2,+}}}_{\text{\shortstack{Previous BESS\\ energy content}}}  - \underbrace{\frac{1}{\eta_{\textrm{B}}} \cdot P_{\textrm{bid}} \cdot \Delta t_{\textrm{2,$\vartheta$}}^{\textrm{+}}}_{\text{\shortstack{If all energy had to\\ be provided by the BESS}}} + \\ 
+ \underbrace{min \left\{ 0.5 \cdot \frac{1}{\eta_{\textrm{B}}} \cdot P_{\textrm{bid}},\sum_{\mathcal{R}} \left( P_{\textrm{j,t+$\vartheta$}}^{\textrm{g,max}} - P_{\textrm{j,t+$\vartheta$}}^{\textrm{g}} \right) \right\}\cdot \Delta t_{\textrm{2,$\vartheta$}}^{\textrm{+}}}_{\text{\shortstack{Part that can be provided by PVs \\(max. 50\% of worst case call)}}} + \\
+ \underbrace{ \eta_{\textrm{B}} \cdot P_{\textrm{j,t+$\vartheta$}}^{\textrm{B,ch}} \cdot \Delta t -\frac{1}{\eta_{\textrm{B}}} \cdot P_{\textrm{j,t+$\vartheta$}}^{\textrm{B,dis}} \cdot \Delta t}_{\text{\shortstack{Scheduled BESS operation}}} 
\end{multline}
\begin{multline}\label{eq:SFC_en_down}
E_{\textrm{\textrm{j,t+$\vartheta$}}}^{\textrm{B,2,-}} = E_{\textrm{j,t+$\vartheta$-1}}^{\textrm{B,2,-}} + \eta_{\textrm{B}} \cdot P_{\textrm{bid}} \cdot \Delta t_{\textrm{2,$\vartheta$}}^{\textrm{-}}\\ 
- min \left\{0.5 \cdot \eta_{\textrm{B}} \cdot P_{\textrm{bid}},\sum_{\mathcal{R}}  P_{\textrm{j,t+$\vartheta$}}^{\textrm{g}} \right\} \cdot \Delta t_{\textrm{2,$\vartheta$}}^{\textrm{-}} \\   + \eta_{\textrm{B}} \cdot P_{\textrm{j,t+$\vartheta$}}^{\textrm{B,ch}}  \cdot \Delta t -\frac{1}{\eta_{\textrm{B}}} \cdot P_{\textrm{j,t+$\vartheta$}}^{\textrm{B,dis}} \cdot \Delta t 
\end{multline}
\noindent where $E_{j,t+\vartheta}^{\textrm{B,2,+}}$ (resp. $E_{j,t+\vartheta}^{\textrm{B,2,-}}$) is the BESS energy content at time t + $\vartheta$ for a call of up (resp. down) regulation at time t; $\vartheta \in \{1,2,3,4\}$ denotes the time for the first 4 hours of the worst case calls, e.g. $E_{j,t}^{\textrm{B,2,+}}$ and $E_{j,t}^{\textrm{B,2,-}}$ correspond to the initial BESS content when the SFC call occurs; and $\Delta t_{\textrm{2,$\vartheta$}}^{\textrm{$\pm$}}$ denotes the worst case up- and down-regulation delivery time at hour $\vartheta$ (whereby $\Delta t_{\textrm{2,$\vartheta$}}^{\textrm{$\pm$}} \leq \Delta t$) derived empirically by the ex-post analysis of the SFC signal~\cite{JannickMA}.

The battery energy content for each individual case, i.e. $\forall t, \vartheta$, are required to stay within the acceptable boundaries,
\begin{subequations}\label{eq:SRL_scenario_bounds_up}
\begin{gather}
SoC_{\textrm{min}}^{\textrm{B}} \cdot E_{\textrm{cap,j}}^{\textrm{B}} \leq E_{\textrm{j,t+$\vartheta$}}^{\textrm{B,2,+}} \leq SoC_{\textrm{max}}^{\textrm{B}} \cdot E_{\textrm{cap,j}}^{\textrm{B}},
\\
\label{eq:SRL_scenario_bounds_down}
SoC_{\textrm{min}}^{\textrm{B}} \cdot E_{\textrm{cap,j}}^{\textrm{B}} \leq E_{\textrm{j,t+$\vartheta$}}^{\textrm{B,2,-}} \leq SoC_{\textrm{max}}^{\textrm{B}} \cdot E_{\textrm{cap,j}}^{\textrm{B}}.
\end{gather}
\end{subequations}

\subsubsection{Tertiary Frequency Control (TFC)}\label{subsubsec:TFC_up}
Tertiary control is asymmetric (up and down) and significantly slower than PFC and SFC, allowing also flexible loads to participate. For this product, both weekly bids as well as bids for single 4-hour blocks can be provided. In the latter case, the constraints apply only to these 4 hours.
The equations are similar to the case of SFC; however, the amount of energy reserves is defined exactly by the regulation of this frequency product, without the need of setting empirical additional constraints. Throughout the duration of the four hours, the full amount of power has to be dispatchable. 

Similar to the case of secondary control, a minimum share of energy has to be provided by the BESS. Here, we define that PV generation combined with flexible loads can account for a maximum share of 80\% of a call. The power and energy constraints for up-regulation are given by
\begin{multline}\label{eq:TFC_power_up}
\sum_{\textrm{j}\in\mathcal{B}} \left( P_{\textrm{j,max}}^{\textrm{B, j}} - P_{\textrm{j,t}}^{\textrm{B,dis}} + P_{\textrm{j,t}}^{\textrm{B,ch}}\right) + \sum_{\textrm{j}\in\mathcal{R}} \left( P_{\textrm{j,t}}^{\textrm{g,max}} - P_{\textrm{j,t}}^{\textrm{g}} \right) + \\
\sum_{\textrm{j}\in\mathcal{F}} \left( f_{\textrm{j,t}}^{\textrm{lflex}} -(-1) \right) \cdot P_{\textrm{j,t}}^{\textrm{shift}} \geq P_{\textrm{bid}},
\end{multline}

\begin{multline}\label{eq:TFC_en_up}
 E_{\textrm{j,t+$\vartheta$}}^{\textrm{B,3,+}} = E_{\textrm{j,t+$\vartheta$-1}}^{\textrm{B,3,+}} - \frac{P_{\textrm{bid}} \cdot \Delta t_{\textrm{3}}}{\eta_{\textrm{B}}} +\\ 
  + min \left\{0.8 \cdot \frac{1}{\eta_{\textrm{B}}} \cdot P_{\textrm{bid}},
\sum_{\textrm{j}\in\mathcal{R}} \left( P_{\textrm{j,t}}^{\textrm{g,max}} - P_{\textrm{j,t}}^{\textrm{g}} \right) +  \right.\\
 \left. + \sum_{\textrm{j}\in\mathcal{F}} \left(f_{\textrm{j,t}}^{\textrm{lflex}} + 1 \right) \cdot P_{\textrm{j,t}}^{\textrm{shift}} \right\}  \cdot \Delta t_{\textrm{3}} + \left( \eta_{\textrm{B}} \cdot P_{\textrm{j,t}}^{\textrm{B,ch}} - \frac{P_{\textrm{j,t}}^{\textrm{B,dis}}}{\eta_{\textrm{B}}} \right) \cdot \Delta t, 
\end{multline}
\noindent where $P_{\textrm{bid}}$ is the weekly or 4-hour block power size of the TFC bid, and $\Delta t_{\textrm{3}}$ is fixed to 1 hour. The SoC constraint $\forall \vartheta \in \{1,2,3,4\}$ is given by
\begin{equation}\label{eq:TRL_scenario_bounds_up}
SoC_{\textrm{min}}^{\textrm{B}} \cdot E_{\textrm{cap,j}}^{\textrm{B}} \leq E_{\textrm{j,t+$\vartheta$}}^{\textrm{B,3,+}} \leq SoC_{\textrm{max}}^{\textrm{B}} \cdot E_{\textrm{cap,j}}^{\textrm{B}}.
\end{equation}
The case of down regulation is similar and straightforward. Finally, for all cases the maximum bid size is constrained by
\begin{equation}\label{eq:bidsize}
0 \leq P_{\textrm{bid}} \leq P_{\textrm{t}}^{\textrm{bid,max}},  
\end{equation}
\noindent where $P_{\textrm{t}}^{\textrm{bid,max}}$ is the maximum power size of the FC product. We use the same variable ($P_{\textrm{bid}}$) for the different FC products, because only one can be offered at a time, i.e. we do not consider provision of multiple services by BESS~\cite{Megel2015}.


\subsection{Islanded mode}\label{subsec:islanded}
In this work, we consider the capability of the distribution grid to be operated in islanded mode, i.e. as a microgrid disconnected from the higher grid level. This is treated by introducing a second set of variables. Most of these constraints are the same as the equations for the grid connected mode and can simply be duplicated.

In this work, the goal in islanded mode is to serve as much of the critical load as possible during the first 24 hours. To achieve that, we utilize the PV generation, BESS and load curtailment.
We treat flexible loads as not critical and thus these loads are not considered in the islanded mode. The power balance equations are given by
\begin{subequations} \label{eq:nodal_active_power_balance_isl}
\begin{gather}
P_{\textrm{j,t}}^{\textrm{inj,isl}} = P_{\textrm{j,t}}^{\textrm{g,isl}} -  \alpha_{\textrm{j,t}}^{\textrm{serv,isl}} \cdot P_{\textrm{j,t}}^{\textrm{l,isl}}  - \left( P_{\textrm{j,t}}^{\textrm{B,ch,isl}} - P_{\textrm{j,t}}^{\textrm{B,dis,isl}} \right) 
\label{eq:nodal_reactive_power_balance_isl},\\
Q_{\textrm{j,t}}^{\textrm{inj,isl}} = Q_{\textrm{j,t}}^{\textrm{g,isl}} + Q_{\textrm{j,t}}^{\textrm{B,isl}} - \alpha_{\textrm{j,t}}^{\textrm{serv,isl}} \cdot P_{\textrm{j,t}}^{\textrm{l,isl}} \cdot \tan(\phi_{l}),\\ \label{eq:nodal_reactive_power_balance_isl}
 0.1 \leq \alpha_{\textrm{j,t}}^{\textrm{serv,isl}} \leq 1, 
 \end{gather}
\end{subequations}
\noindent where $\alpha_{\textrm{j,t}}^{\textrm{serv,isl}}$ denotes here the fraction of active power served. 

Modern grid codes require a minimum power factor requirement in the grid-connected case\cite{VDE}. However, in the islanded mode we exploit the full functionality of the PV and BESS inverters. Thus, the reactive power provision is described by
\begin{equation}\label{eq:react_generation_isl}
(Q_{\textrm{j,t}}^{\textrm{g,isl}})^{2} \leq (S_{\textrm{j,t}}^{\textrm{g,isl}})^{2} - (P_{\textrm{j,t}}^{\textrm{g,isl}})^{2}.
\end{equation}

Finally, all the constraints concerning the OLTC are not active in the islanded case. The only link between the set of variables in the grid-connected and the islanded mode is the BESS energy content at timestep $\tau$, when the islanding operation begins, i.e.$E_{\textrm{j,$\tau$}}^{\textrm{B,isl}}=E_{\textrm{j,$\tau$}}^{\textrm{B}}$. After that, the two sets of variables describe independent possible future developments.
\review{\section{Handling of uncertainty and chance-constrained OPF formulation}\label{Sec:uncertainty_finalproblem}
This section first describes how the uncertainty is considered in form of chance constraints and then summarizes the final centralized CC-OPF formulation.}
\subsection{Accounting for Uncertainty through Chance Constraints}\label{uncertain}
In order to consider the impact of generation uncertainty, we follow our previous work~\cite{StavrosIREP,stavrosML} and we re-formulate the problem using chance constraints~\cite{roald2017submitted,Roald2013}. We assume that the PV power injection is the only source of uncertainty (load uncertainty can be also included in a similar way) and we use as input forecast error distributions with different forecasting horizons (1 to 24 hours ahead).

Following~\cite{StavrosIREP,stavrosML} we model the voltage and current constraints as chance constraints that will hold with a chosen probability $1-\varepsilon$, where $\varepsilon$ is the acceptable violation probability. E.g., the voltage and current magnitude constraints are reformulated as $\mathbb{P}\left\{ V_{\textrm{min}} \leq |V_{\textrm{j,t}}| \leq  V_{\textrm{max}} \right\}\geq{1-\varepsilon}$ and  $\mathbb{P}\left\{|I_{\textrm{i,t}}^{\textrm{br}}| \leq  I_{\textrm{i}^{\textrm{max}}} \right\}\geq{1-\varepsilon}$, respectively.   
To solve the resulting CC-OPF, we interpret the probabilistic constraints as \emph{tightened} deterministic versions of the original constraints following the work of~\cite{Roald2013,roald2017submitted}. The tightening represents a security margin against uncertainty, i.e., an \emph{uncertainty margin}. Thus, we express the voltage and current constraints as
\begin{align}
    \begin{cases}
               |V_{\textrm{j,t}} | & \leq V_{\textrm{max}} - \Omega_\textrm{V j,t}^\textrm{upper}\\
               {Re}\left\{V_{\textrm{j,t}} \right\} &\geq V_\textrm{min} + \Omega_\textrm{V j,t}^\textrm{lower}, 
    \end{cases} \label{eq:vol_lim2}\\
    |I_\textrm{br,i,t}|  \leq I_{\textrm{i,max}} - \Omega_{I_\textrm{br,i}}, \label{eq:current_lim2}
\end{align}
where $\Omega_\textrm{V}^\textrm{lower},~\Omega_\textrm{V}^\textrm{upper}$ are the tightenings for the lower and upper voltage magnitude constraints and $\Omega_{I_\textrm{br}}$ are the tightenings of the current magnitude constraints. The interested reader is referred to~\cite{StavrosIREP} for more details on this method.

The uncertainty margins are constant within the OPF solution process, and evaluated outside of the OPF solution. Thus, we use a Monte Carlo approach and the non-linear AC power flow equations to evaluate the boundaries. This further allows us to include any uncertainty probability distribution. 

Hence, we form empirical distributions for the voltage and current chance constraints at each time step based on the results from the Monte Carlo simulations. To enforce a chance constraint with $1-\epsilon$ probability we need to ensure that the $1-\epsilon$ quantile of the distribution remains within the bounds. Thus, the tightening corresponds to the difference between the forecasted value with zero forecast error and the $1-\epsilon$ quantile value evaluated based on the empirical distribution resulting from the Monte Carlo Simulations, e.g. $|V_{\textrm{bus,j,t}}^{\textrm{0}}|$ and $|V_{\textrm{bus,j,t}}^{\textrm{1-$\epsilon$\%}}|$ for the voltage constraints.
The empirical uncertainty margins to be used in the next iteration are then given by 
\begin{subequations}
\begin{align}
    \Omega_\textrm{V j,t}^\textrm{upper}   &= |V_{\textrm{bus,j,t}}^{\textrm{1-$\epsilon$}}| - |V_{\textrm{bus,j,t}}^{\textrm{0}}|,\\
    \Omega_\textrm{V j,t}^\textrm{lower}   &= |V_{\textrm{bus,j,t}}^{\textrm{0}}|  - |V_{\textrm{bus,j,t}}^{\textrm{$\epsilon$}}|,\\
    \Omega_{I_\textrm{br,i}}^\textrm{upper}&= |I_\textrm{br,i,t}^{\textrm{1-$\epsilon$}}| - |I_\textrm{br,i,t}^{\textrm{0}}|,
\end{align}
\end{subequations}
where superscript $^0$ indicates the current or voltage magnitude at the operating point with zero forecast error.
Finally, an iterative algorithm is needed, because the uncertainty margins rely on the derived DER setpoints~\cite{Schmidli2016, roald2017submitted}. Consequently, we alternate between solving a deterministic OPF with tightened constraints, and calculating the uncertainty margins $\Omega_\textrm{V}^\textrm{lower},~\Omega_\textrm{V}^\textrm{upper},~\Omega_{I_\textrm{br}}^\textrm{upper}$. When the change in the tightening values between two subsequent iterations is below a threshold $(\eta_V^{\Omega},~\eta_I^{\Omega})$, then the algorithm has converged.

\subsection{Solution Algorithm}
\label{SolAlgo}
\begin{figure}[]
\centering
\includegraphics[width=0.99\columnwidth]{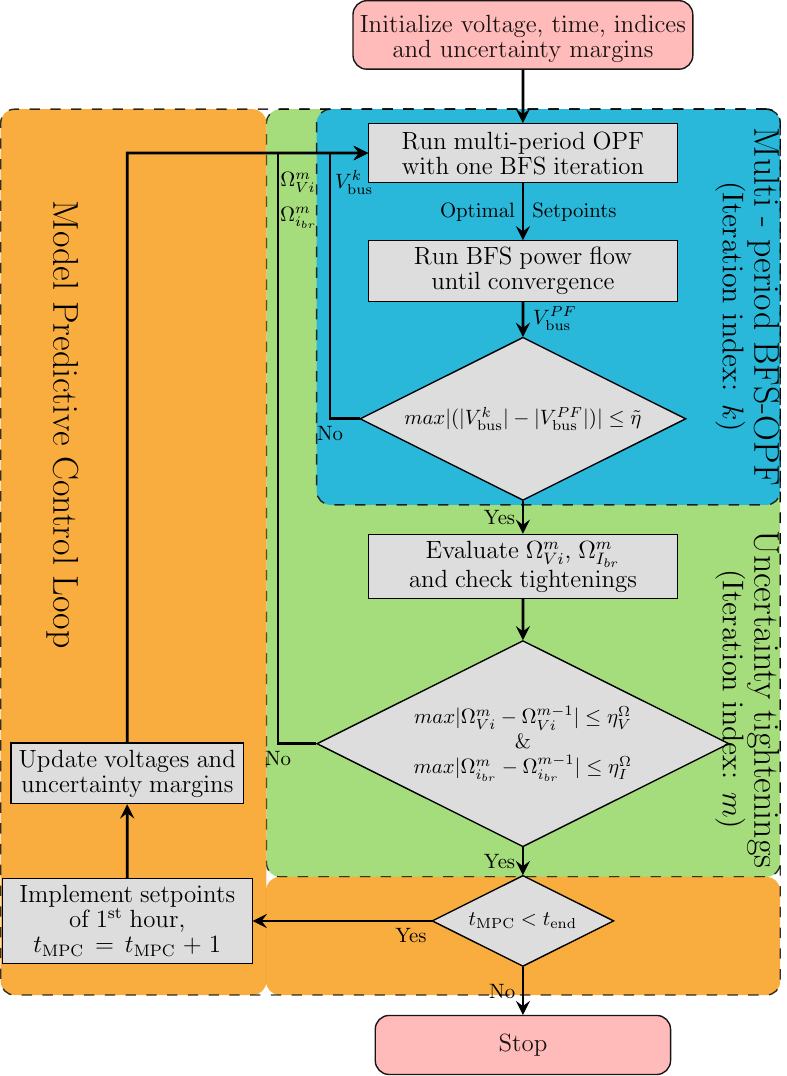}
\caption{Proposed CC-OPF implemented in an MPC fashion} 
\label{fig:SolAlgo}
\end{figure}

In this section, we summarize the proposed solution method for the centralized CC-OPF scheme implemented in an MPC fashion, sketched in Fig.~\ref{fig:SolAlgo}. First, the initialization stage sets the uncertainty margins to zero and initializes the voltage levels to a flat voltage profile. At the core of the proposed methodology lies the formulation of the multi-period centralized CC-OPF, which considers the provision of ancillary services as well as the possibility for islanded operation.
The CC-OPF calculates the optimal DER setpoints based on a single sweep of the BFS algorithm. The BFS power-flow algorithm then runs until convergence for the obtained control settings. The CC-OPF is then performed again using the updated voltages from the full BFS. These inner iterations are carried out until convergence. After the multi-period OPF has converged, the uncertainty margins are evaluated in the outer loop as described in Section~\ref{uncertain}. The iteration index of the OPF loop is denoted by $k$ and the iteration of the uncertainty loop by $m$. The iterative procedure continues until all parts of the algorithm have reached convergence. Then, only the optimal setpoints of the first
time step are implemented. Subsequently, the PV forecast is updated,
the current timestep is increased and the next CC-OPF problem with a horizon of 24 hours is solved.

\review{The resulting optimization problem is a mixed-integer quadratically constrained program (MIQCP) and can be solved efficiently by modern powerful solvers. The computational burden depends on the dimensions of the grid, the acceptable violation probability, and the number and complexity of the considered DGs.}

\review{Due to the efficient handling of the power flow equations through the BFS formulation, hundreds of nodes and branches can be handled without a drastic increase in the computational burden. Regarding the uncertainty handling, the selection of $\epsilon$ influences the execution time of the proposed scheme, since it modifies the feasible area of the optimization problem. The larger the required fulfillment (small values of epsilon), the smaller the feasible area of the optimization problem, making the optimization more demanding. If it is necessary to reduce the computational burden, DGs with complex modes can be handled with reasonable approximations. E.g. constraint~\eqref{eq:BESS_P2} could be replaced as in~\cite{stavrosPowertech} to avoid the need for binary variables, and the operation of the tap changers could be modeled as continuous variables, rounded ex post to the closest integer.}

\review{Overall, however, realistic distribution grid dimensions require solving time in the range of minutes, which is acceptable for such kind of steady state analysis and can be implemented in existing active distribution grids.}

\review{}

\section{Case Study - Results}\label{case}
\begin{figure}[]
    \begin{centering}
	\includegraphics[width=0.8\columnwidth]{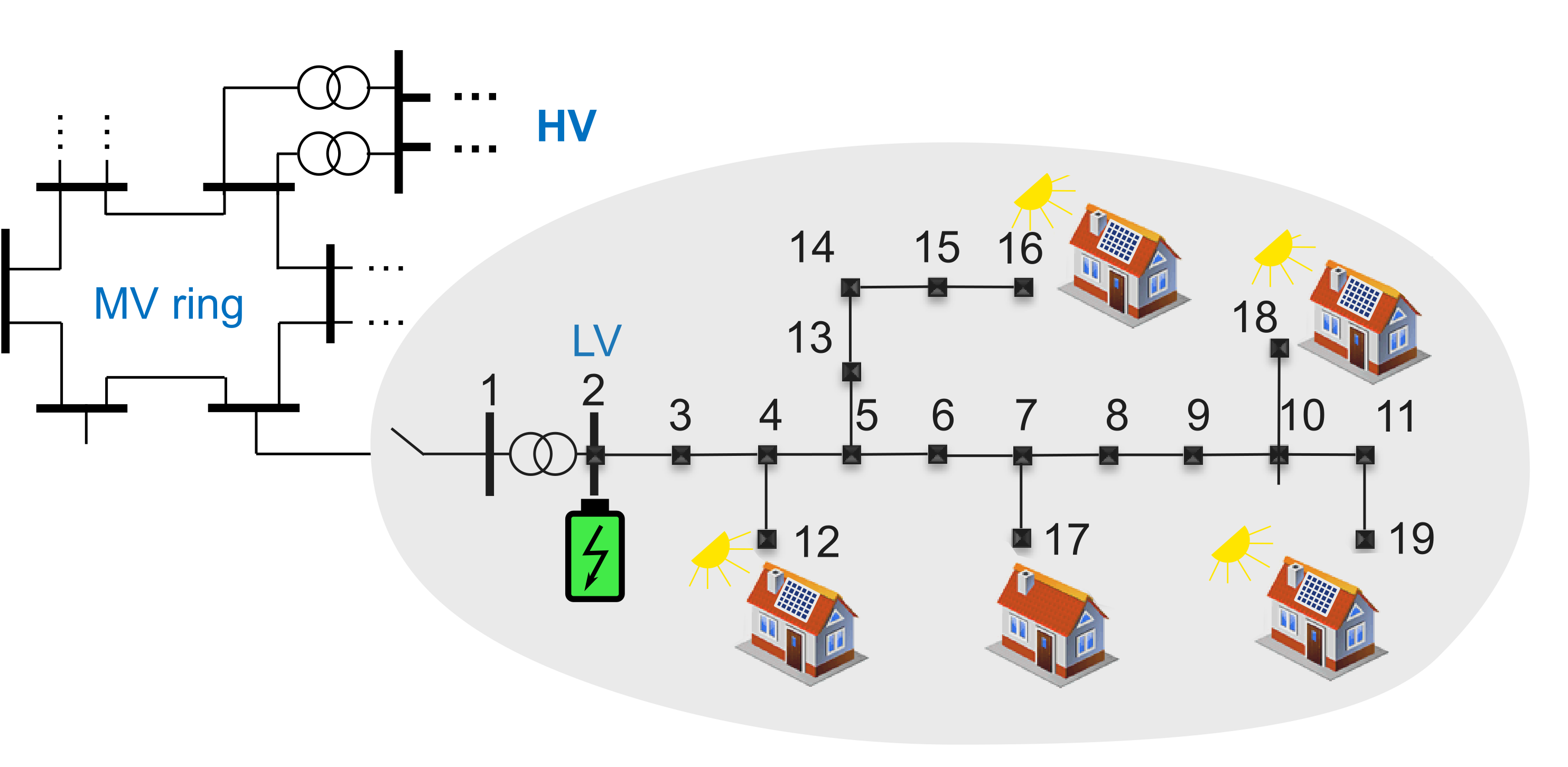}
	 \caption{Cigre residential European LV grid able to operate in grid-connected and islanded mode.}	\label{fig:cigre_test_system} 
	\end{centering}
	\vspace{-0.2cm}
\end{figure}

In order to demonstrate the proposed method, we use a typical European radial LV grid~\cite{Strunz2014}, sketched in Fig.~\ref{fig:cigre_test_system}. The installed PV capacity is expressed as a percentage of the total peak load as follows: PV nodes = [12, 16, 18, 19], PV share~(\%) = [35, 25, 30, 45]. Furthermore,
we consider flexible loads up to 5 kW at nodes [17, 18, 19], \review{i.e. $5\%$, $15\%$ and $10\%$ of the corresponding nominal load}. The BESS capacity at node 2 is 484 kWh, \review{and the maximum power 484 kW}. In this work, we only consider balanced, single-phase system operation, but the framework
can be extended to three-phase unbalanced networks as we explain in~\cite{StavrosPSCC18}.

The spot market prices were assumed  equal to the realized values of 2016~\cite{EPEX}. The realized reserve prices of 2016 are available in~\cite{swissgrid}. \review{To adjust the cost for the islanded case, we used a constant of $f_{\textrm{isl}}=0.1$, and very high load curtailment cost of $c_{\textrm{t}}^{\textrm{curt,l,isl}}=250\frac{\textrm{\euro}}{\textrm{MWh}}$.} Furthermore, a realized primary control signal was derived from a frequency signal with a temporal resolution of one second. A realized secondary control signal with the same time resolution was taken from~\cite{ekz}.

Regarding the uncertainty modeling, we use historical forecast
error distributions from an area in Switzerland provided by~\cite{ekz} and we enforce the chance constraints with an $\epsilon= 5\%$ violation probability. We assume a maximum
acceptable voltage of $1.1$ p.u and cable current magnitude of
$1$ p.u. on the cable base. The minimum acceptable voltage
is set to $0.9$ p.u..

Using this system, we investigate the capability of the DN to switch to islanded mode, while offering frequency control products. Furthermore, we show how the DN responds to a frequency control call, respecting the islanding requirement.
The implementation was done in MATLAB. For the centralized OPF-based control, YALMIP~\cite{Lofberg2004} was used as the modeling layer and Gurobi~\cite{gurobi} as the solver. The results were obtained on an Intel Core i7-2600 CPU and 16 GB of RAM.

\subsection{Islanded operation}\label{sec:isl_operation}
The first part of the results refers to the ability of the DN to switch to the islanded mode, where at least 10$\%$ of the load should be served for the next 24 hours. This parameter is estimated to cover emergency services.\\
\subsubsection{Determination of minimal BESS size}\label{sec:min_BESS_size}
\begin{figure}[]
\centering
\includegraphics[trim=1.5cm 8.5cm 0.8cm 9.0cm, clip=true, width=0.9\columnwidth]{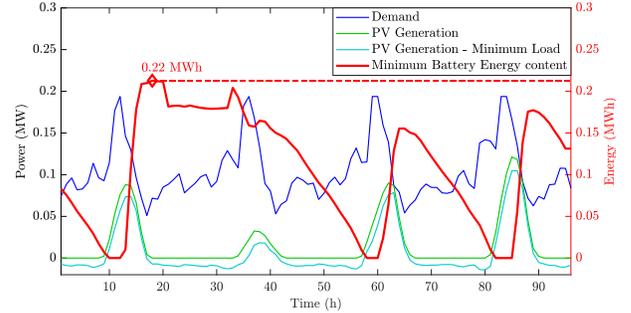}
\caption{Historical worst case conditions to determine the minimum BESS energy capacity}
\label{fig:wc_BESS}
\vspace{-0.4cm}
\end{figure}

A minimum battery energy capacity is required in order to ensure islanded feasibility under different PV injection and loading conditions. Thus, we used historical values of available PV and load data to determine the minimum BESS requirement for islanded operation. We performed yearly MPC-OPF calculations with a 24-hour horizon, without considering uncertainties, to estimate the needed BESS size iteratively; i.e. we kept increasing the BESS size until we derived feasible solutions for the whole year. The worst case period is shown in Figure \ref{fig:wc_BESS}, indicating a minimum BESS size of $220$ kWh. 

In order to allow provision of AS, we investigated various BESS capacities corresponding to $1.4-2.6$ times the needed minimum value. In the remaining simulations, we will consider a BESS of $484$ kWh.
\subsubsection{Switch to islanded mode}\label{subsec:balancing_isl}
\begin{figure}[]
\centering
\includegraphics[trim=3.3cm 10.3cm 3.3cm 10.0cm, clip=true, width=0.95\columnwidth]{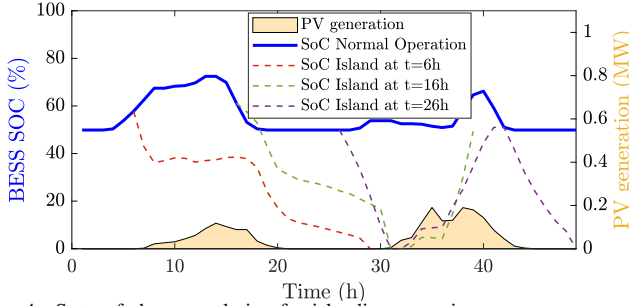}
\caption{State of charge evolution for islanding scenarios}
\label{fig:isl_feas}
\vspace{-0.2cm}
\end{figure}

According to Section~\ref{sec:isl_operation}, the switch to islanded operation should be feasible at any time instant. Figure \ref{fig:isl_feas} shows the evolution of the BESS SOC for islanding at distinct hours in the considered time period. The power balance is kept using the BESS capacity, PV injections, and load and PV curtailment. As can be observed, the BESS SoC evolution depends on the PV generation and load forecasts; At noon hours, the PV units provide power for the loads and BESS charging, while at night the BESS is discharged to guarantee a 24-hour islanded operation. 
\vspace{-0.2cm}
\subsection{Frequency control}
In the grid-connected case, the DN offers frequency regulation as an ancillary service, while at the same time fulfilling the islanded requirement for the next 24 hours. 

\subsubsection{PFC}

\review{Assuming that the BESS is always charging or discharging at a maximum rate of 1C\footnote{A C-rate is a measure of the rate at which a BESS is charged or discharged relative to its maximum capacity. A 1C rate means that the discharge current will discharge the entire
battery in 1 hour.} to limit the capacity fading from offering frequency control products~\cite{StavrosPSCC2016}, an energy requirement of 15 minutes PFC power in both directions, i.e. 30 minutes in total, translates into reserving 50\% of the total BESS storage capacity.} The algorithm keeps the SoC at the upper limit to minimize load shedding in case of a switch to the islanded mode.

\begin{figure}[]
\centering
\includegraphics[trim=3.5cm 10.0cm 4cm 10.5cm, clip=true, width=0.9\columnwidth]{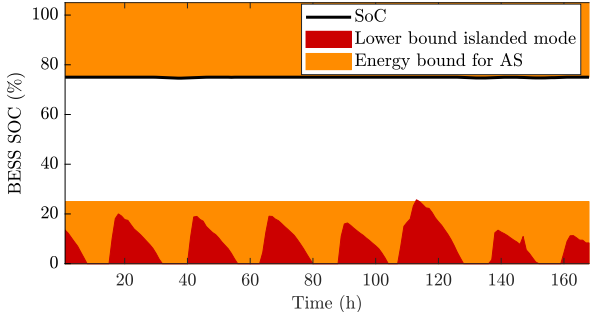}
\caption{BESS SoC with PFC reserve provision}
\label{fig:PFC}
\end{figure}

Figure~\ref{fig:PFC} shows the BESS SoC while providing PFC over a summer week. Staying outside of the red area guarantees that in the case of a switch to islanded operation at any time step the critical load can be supplied by preserving a minimum BESS energy content based on load and PV generation forecasts. The orange area represents the energy limit imposed by the offered frequency control product. The white area defines the allowable feasible region for the SOC, with the black line showing the optimization result. In case of overlapping between the orange and red area, the more limiting area is relevant. In case of operating in islanded mode, frequency reserves are not provided anymore.

\subsubsection{SFC}
\begin{figure}[]
\centering
\includegraphics[trim=3.5cm 10.0cm 4cm 10.5cm, clip=true, width=0.9\columnwidth]{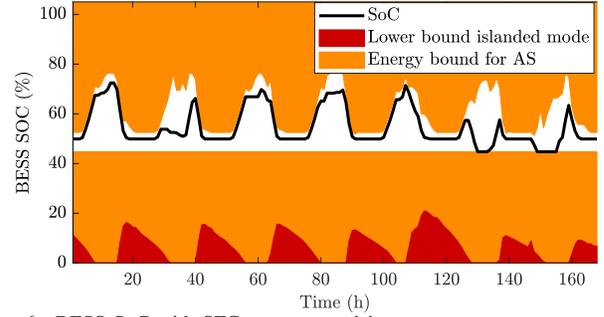}
\caption{BESS SoC with SFC reserve provision}
\vspace{-0.2cm}
\label{fig:SFC}
\end{figure}

For this product we consider also PV units, which can curtail power providing down-regulation. Hence, the upper bound on the energy level of the storage during hours with PV injections is relaxed, as seen in Fig.~\ref{fig:SFC}. The BESS can be charged during these hours, leading to higher self-consumption and more available energy in case of a switch to islanded mode.

\subsubsection{TFC - weekly offer} 
\begin{figure}[]
\vspace{-0.3cm}
\centering
\includegraphics[trim=3.5cm 10.0cm 4cm 10.5cm, clip=true, width=0.9\columnwidth]{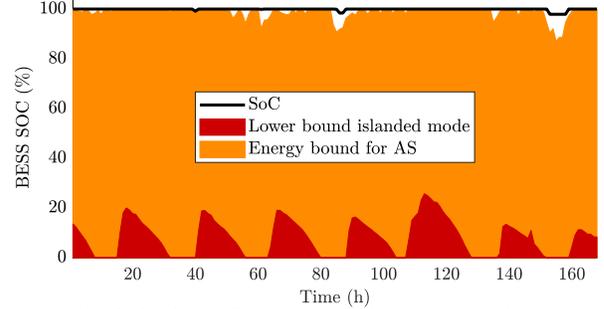}
\caption{BESS SoC with TFC reserve provision - up regulation}
\label{fig:TFC}
\end{figure}

\paragraph{Up regulation} Providing maximum up TFC regulation resulted in a fully charged BESS, as can be observed in Fig.~\ref{fig:TFC}. In this way, we not only achieve maximum reserve provision, but also minimum load curtailment in the islanded mode. Limited flexibility is offered by flexible loads, as can be seen by the white areas, the size of which does not influence the maximum bid size.

\begin{figure}[]
\vspace{-0.5cm}
\centering
\includegraphics[trim=3.5cm 10.0cm 4cm 10.5cm, clip=true, width=0.9\columnwidth]{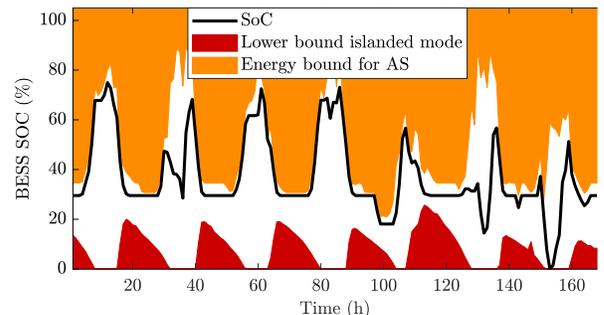}
\caption{BESS SoC with TFC reserve provision - down regulation}
\label{fig:TFCdown}
\vspace{-0.3cm}
\end{figure}

\paragraph{Down regulation}
The case of down regulation is shown in Fig.~\ref{fig:TFCdown}, where the optimization tries to keep the SoC low in order to respond to a TFC dispatch call, while at the same time respecting the islanding requirement. Similar to the SFC case, during noon hours with solar power, the SoC can be increased, since PV power curtailment is available.

\subsection{Call for SFC}\label{subsubsec:SRL_call}
So far, we studied the needed power and energy reserves. In this section, we simulate the response of the DN to an actual continuous SFC signal. Since we cannot forecast the signal, we used the realized signal from 2016. 
Figure~\ref{fig:srl_call_E} shows the worst-case week in terms of needed power of the SFC signal as well as the corresponding cumulative energy requirement. We consider the possibility of participating in the spot market with a lead time of four hours. As can be observed, the algorithm chooses to buy energy on the spot market three times indicated by the red circles. 
The dashed line corresponds to the cumulative energy demand without spot market purchases, whereas the solid line to the resulting cumulative energy given the purchases in the spot market.
\review{Finally, Fig.~\ref{fig:srl_call_SOC} shows the evolution of the BESS SoC following the SFC signal in solid, and the SoC without offering SFC with a dashed line. As can be observed, the three purchases of power are needed so that the SoC is kept high enough to allow for the islanded mode.
\begin{figure}[]
\centering
\includegraphics[width=1\columnwidth]{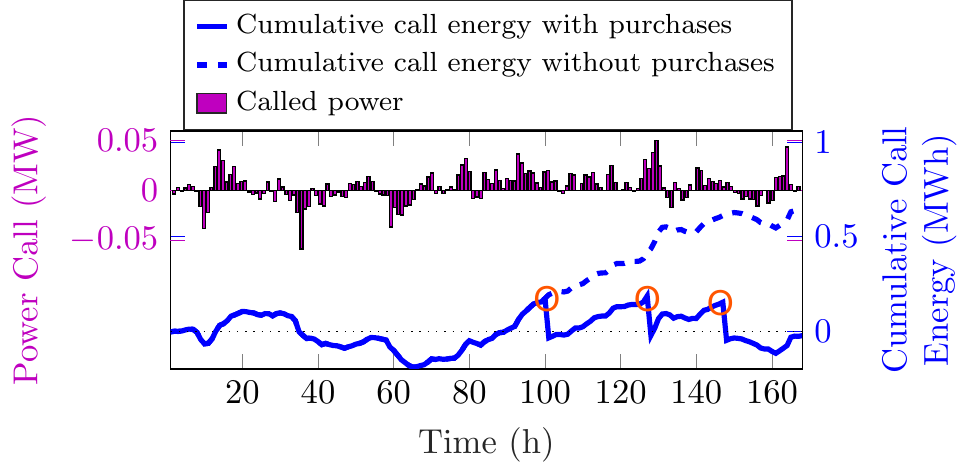}
\caption{Secondary control call signal and cumulative energy requirement for $484 \textrm{kWh}$ BESS in Summer}
\label{fig:srl_call_E}
\end{figure}
 \begin{figure}[]
 \vspace{-0.2cm}
 \centering
 \includegraphics[width=1\columnwidth]{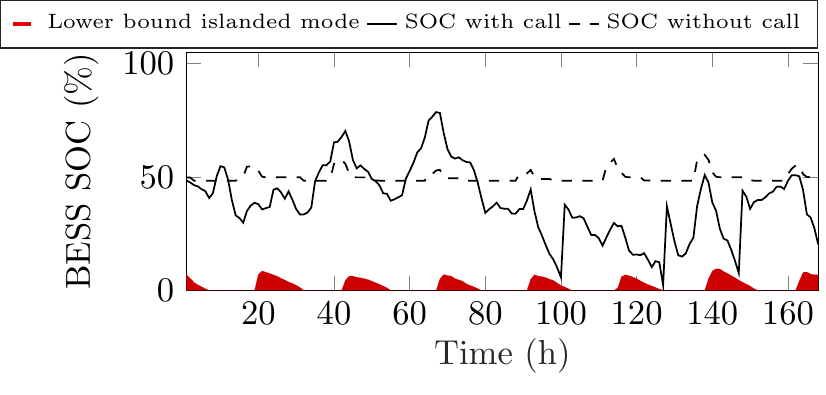}
 \caption{BESS SoC with secondary call signal}
 \label{fig:srl_call_SOC}
 \vspace{-0.3cm}
 \end{figure}
}
\review{\subsection{Impact of BESS size on the rating of the transformer}\label{subsec:MinTrafo}
As a final case study, we investigate the impact of the BESS size on the needed rating of the MV/LV transformer, without offering frequency control products. The BESS can contribute to the power needed to and from the active distribution grid, reducing the required transformer rating. In this way, the service of investment deferral can be offered to the operator, which might need to cope with increasing demand or DG injections.} 

\begin{figure}
\vspace{-0.5cm}
\centering
\includegraphics[trim=3.5cm 9.7cm 5cm 9.5cm, clip=true, width=1\columnwidth]{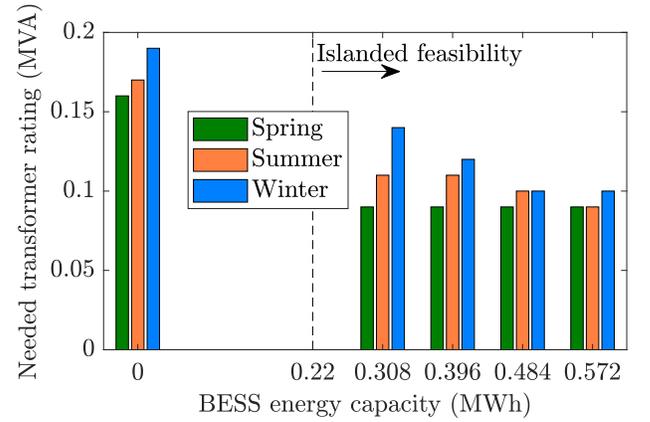}
\caption{BESS SoC with TFC reserve provision - down regulation}
\label{fig:Min_Trafo}
\vspace{-0.5cm}
\end{figure}

\review{Figure~\ref{fig:Min_Trafo} shows the required rating of the secondary substation transformer, varying the energy capacity of the BESS placed at the same node. A seasonal analysis allows calculating the most critical period, i.e. winter in our case, that defines the needed rating. We observe that the larger the energy BESS capacity, the smaller the required transformer rating; however, the BESS contribution is decreasing with increasing BESS size.} 

\section{Conclusion}  \label{Conclusion}
Modern DNs consider the active control capabilities of DERs in order to provide a secure, reliable and optimal operation of the grid. Furthermore, they can offer ancillary services to higher voltage levels, or even operate disconnected from the main grid. 

In this paper, we have shown that ADGs can be coordinated through centralized control schemes to provide ancillary services and provision for islanded operation. The proposed method allows ADGs to support the transmission network but at the same time provide increased resilience through controlled islanding. We have shown how the different operational requirements can be formulated in the problem constraints and provided techniques to tackle the uncertainty.

\section*{Acknowledgments}
\review{The work of P. Aristidou was partially supported by the Engineering and Physical Sciences Research Council (EPSRC) in the UK under grant reference EP/R030243/1.}
\bibliographystyle{IEEEtran}
\bibliography{bibliography}

\end{document}